\documentstyle[12pt,amssymb]{article}

   \makeatletter
   \def\maketitle{\begin{center}\let\thanks=\footnote
      {\large\bf\@title}\par\bigskip\bigskip{\sc\@author}\par\bigskip{\rm\@date}
      \end{center}\bigskip\thispagestyle{empty}}

   \def\@begintheorem#1#2{\trivlist\item[\hskip\labelsep{\bf#1 #2.}]\it}
   \def\@opargbegintheorem#1#2#3{%
      \trivlist\item[\hskip\labelsep{\bf#1 #2 {\rm(#3)}.}]\it}

   \let\o@item=\@item
   \def\@item[#1]{\o@item[\rm #1]}

   \def\@biblabel#1{#1.}
   \let\othebibliography=\thebibliography
   \def\thebibliography#1{\small
      \def\@listi{\topsep=0cm\parsep=0cm\itemsep=0cm}\othebibliography{#1}}

   \newtheorem{satz}{Satz}[section]
   \makeatletter
   \@addtoreset{equation}{satz}
   
   \newtheorem{theorem}[satz]{Theorem}
   \newtheorem{lemma}[satz]{Lemma}
   \newtheorem{proposition}[satz]{Proposition}
   \newenvironment{varthm*}[1]{\trivlist\item[]\it{\bf #1.}}{\endtrivlist}
   \def\startproof{\addvspace{\bigskipamount}\noindent}
   \def\proof{\startproof{\it Proof. }}
   
   \def\qedsymbol{\frame{\rule[0pt]{0pt}{8pt}\rule[0pt]{8pt}{0pt}}}
   \def\qed{\nopagebreak\hspace*{\fill}\qedsymbol\par\addvspace{\bigskipamount}}

   \def\to{\longrightarrow}
   \def\mapsto{\mapstochar\longrightarrow}
   \def\epsilon{\varepsilon}
   \def\hat{\widehat}
   \def\bar{\overline}
   \def\({\left(}
   \def\){\right)}

   \def\subjclass#1{{\def\thefootnote{}
      \footnotetext{1991 {\it Mathematics Subject Classification:} #1.}}}

   \newenvironment{items}{\list{$\bullet$}{
      \parsep=0cm\itemsep=0cm\topsep=0cm\partopsep=\medskipamount
      \def\makelabel##1{\hss\llap{##1}}}
   }{\endlist}

   \def\bbQ{{\Bbb Q}}
   \def\bbR{{\Bbb R}}
   \def\bbC{{\Bbb C}}
   \def\bbZ{{\Bbb Z}}

   \def\be{\begin{eqnarray*}}
   \def\ee{\end{eqnarray*}}
   
   \def\Bigwith{\ \Big\vert\ }
   \def\liste#1#2#3{\mbox{$#1_{#2},\dots,#1_{#3}$}}
   \def\eqnref#1{(\ref{#1})}
   
   \def\tensor{\otimes}
   
   \def\with{\mid}
   \def\isom{\cong}
   \def\inverse{^{-1}}
   \def\operatorname#1{\mathop{\rm #1}\nolimits}
   \def\im{\operatorname{im}}
   \def\End{\operatorname{End}}
   \def\EndQ{\End_{\bbQ}}
   \def\NSQ{\NS_{\bbQ}}
   \def\Pic{\operatorname{Pic}}
   \def\rk{\operatorname{rk}}
   \def\eqdef{=_{\operatorname{def}}}

   \def\Nef{\operatorname{Nef}}
   \def\N{N}
   \def\NE{N\hskip-0.25em E}
   \def\bNE{{}\hskip0.25em\bar{\hskip-0.25em N\hskip-0.25em E}}
   \def\NS{N\hskip-0.2em S}
   \def\NSR{\NS_{\bbR}}
   \def\NP{\NS^+}
   \def\NX{\NP(X)}
   \def\and{\quad\mbox{ and }\quad}
   \def\implies{\Rightarrow}
   \def\i{\iota_1^*}
   \def\ii{\iota_2^*}
   \def\iii{\iota_3^*}

\begin{document}

\title{On the cone of curves of an abelian variety}
\author{Thomas Bauer}
\date{}
\maketitle


\setcounter{section}{-1}
\section{Introduction}

   Let $X$ be a smooth projective variety over the 
   complex numbers and let $\N_1(X)$ be the real vector
   space
   $$
      \N_1(X)\eqdef\{\mbox{1-cycles on $X$ modulo numerical
      equivalence}\}\tensor\bbR \ .
   $$
   As usual denote by $\NE(X)$ the {\em cone of curves} on $X$,
   i.e.\ the convex cone in
   $\N_1(X)$ generated by the effective 1-cycles.
   The {\em closed cone of curves}
   $\bNE(X)$ is the closure of $\NE(X)$ in $N_1(X)$.
   One knows by the Cone Theorem \cite{Mor82} that
   it is
   rational polyhedral whenever $c_1(X)$ is ample.  For
   varieties $X$ such that
   $c_1(X)$ is not ample, however,
   it is in general difficult to determine the
   structure of $\bNE(X)$, since it may depend in a subtle way
   on the geometry of $X$ (cf.\ \cite[\S4]{CKM}).
   This becomes already apparent in the
   surface case, as work of Kov\'acs on K3 surfaces shows
   (see \cite{Kov94}).

   The purpose of this paper is to study
   the cone of curves of abelian
   varieties.  Specifically, we focus on the problem of
   determining the abelian varieties $X$ 
   such that the closed cone $\bNE(X)$ is rational polyhedral.
   Attacking the question from the dual point of view, one
   is lead to
   consider the nef cone of $\Nef(X)$ or the 
   semi-group $\NX$ of homology classes of
   effective line bundles, i.e.\ the subset
   $$
      \NX\eqdef\{\lambda\in\NS(X)\with\lambda=c_1(L)\mbox{ for some
      }L\in\Pic(X)\mbox{ with }h^0(X,L)> 0\}
   $$
   of the N\'eron-Severi group of $X$.  
   In fact, Rosoff has studied this semi-group in
   \cite{Ros81}, where he gives examples  
   of abelian varieties for which $\NX$ is finitely generated,
   as well as examples where finite generation fails.  
   He shows:
\begin{items}
   \item[(1)] 
      If $X$ is a singular abelian variety, i.e.\ if
      $\rk\NS(X)=(\dim X)^2$, and if $\dim X\ge 2$, 
      then $\NX$ is not finitely generated.
   \item[(2)]
      For elliptic curves $E_1$ and $E_2$, $\NP(E_1\times
      E_2)$ is finitely generated if and only if
      $\rk\NS(E_1\times E_2)=2$.
\end{items}

   Considering these examples it is natural to 
   ask if the abelian varieties $X$ for which $\NX$ is
   finitely generated or, equivalently, for which 
   $\bNE(X)$ is rational polyhedral,
   can be characterized in a simple way.
   Our main result shows that this is in fact the case:

\pagebreak
\begin{varthm*}{Theorem}
   Let $X$ be an abelian variety over the field of complex
   numbers.  Then the following conditions are equivalent:
   \begin{items}
      \item[(ia)]
	 The closed cone of curves $\bNE(X)$ is rational polyhedral.
      \item[(ib)]
	 The nef cone $\Nef(X)$
         is rational polyhedral.
      \item[(ic)]
	 The semi-group $\NX$ is finitely generated.
      \item[(ii)]
	 $X$ is isogenous to a product
	 $$
	    X_1\times\dots\times X_r
         $$
	 of mutually non-isogenous abelian varieties $X_i$
	 with $\NS(X_i)\isom\bbZ$ for $1\le i\le r$.
   \end{items}
\end{varthm*}

   Note that, since on abelian varieties the nef cone
   coincides with the effective cone, the equivalence of 
   (ia), (ib) and (ic) follows from elementary properties of cones
   and is stated here merely for the sake of completeness (see
   Sect.\ \ref{sect proof}).
   Observe that the theorem of course contains statement (2)
   above, while statement (1) follows from the theorem plus
   the fact that by \cite{ShiMit74}
   a singular abelian variety is isogenous to a
   product $E^n$ for some elliptic curve $E$.

\begin{varthm*}{\it Notation and Conventions}
   \rm
   We work throughout over the field $\bbC$ of complex numbers.

   We will always
   use additive notation for the tensor product of line
   bundles, since this is more convenient for our purposes
   (for example when working with $\bbQ$- or $\bbR$-line
   bundles).
   Numerical equivalence of divisors or line bundles, which
   for abelian varieties coincides with algebraic equivalence,
   will be denoted by $\equiv$.

   If $X$ is an abelian variety and $L$ is a line bundle on
   $X$, then $\phi_L$ denotes the homomorphism $X\to\hat X$,
   $x\mapsto t_x^*L - L$, where $t_x$ is the translation map
   $y\mapsto x+y$ and $\hat X=\Pic^0(X)$ is the dual abelian
   variety.  Recall that $\phi_L$ depends only on the
   algebraic equivalence class of $L$.
\end{varthm*}


\section{Effective classes on simple abelian varieties}

   In this section we consider the semi-group of
   effective divisor classes on
   simple abelian varieties.  We start by stating alternative
   characterizations of $\NX$ which we will use in the sequel.
   While these are at least implicitly well-known, we include
   a proof for the convenience of the reader.

\begin{lemma}\label{nef and effective}
   Let $X$ be an abelian variety of dimension $n$ and let $A$
   be an ample line bundle on $X$.  Then the following
   conditions on a line bundle $L\in\Pic(X)$ are equivalent:
   \begin{items}
      \item[(i)] 
	 $L$ is algebraically equivalent to some effective 
	 line bundle.
      \item[(ii)]
	 $L$ is nef.
      \item[(iii)]
	 $L^iA^{n-i}\ge 0$ for $1\le i\le n$.
   \end{items}
\end{lemma}

\proof
   Condition (i) certainly implies (ii), since if $L\equiv L'$
   for some effective line bundle $L'$, then a suitable
   translate of an effective divisor in $|L'|$ will intersect
   any given curve properly.  
   The implication (ii) $\implies$ (iii) is clear, since an
   intersection product of nef line bundles is non-negative.
   For (iii) $\implies$ (ii) it is enough to show that the
   line bundle $A+mL$ is ample for all $m\ge 0$.  But this
   follows from
   $$
      (A+mL)^i A^{n-i}=A^n+\sum_{k=0}^{i-1} {i \choose k}
      A^{n+k-i}L^{i-k} m^{i-k} > 0 \ ,
   $$
   and the version of the Nakai-Moishezon Criterion given in
   \cite[Corollary 4.3.3]{LB}.
   Finally, for the implication (ii) $\implies$ (i), suppose
   that $L$ is nef.  Then $A+mL$ is ample for all $m\ge 0$, so
   that the first Chern class of $L$, viewed as a hermitian
   form on $T_0X$, cannot have negative
   eigenvalues.  But this implies that there is a line bundle
   $P\in\Pic^0(X)$ such that $L+P$ descends to an ample line
   bundle on a quotient of $X$ and is therefore effective
   (cf.\ \cite[Sect.\ 3.3]{LB} and \cite[p.\ 95]{Mum70}).
\qed

   We show next that on a simple abelian variety the existence
   of two algebraically independent line bundles already
   prevents $\NX$ from being finitely generated:

\begin{proposition}\label{simple case}
   Let $X$ be a simple abelian variety such that $\NX$ is
   finitely generated.  Then $\NS(X)\isom\bbZ$.
\end{proposition}

\proof
   Assume to the contrary that $\rk\NS(X)>1$ and choose ample
   line
   bundles $L_1$ and $L_2$ whose classes are not proportional
   in $\NSQ(X)$.  Consider then the positive real number
   $$
      s \eqdef\inf\Big\{t\in\bbR\Bigwith tL_1-L_2\mbox{ is nef }\Big\} \ .
   $$
   Here $tL_1-L_2$ is considered as an $\bbR$-line bundle and
   nefness means that $tL_1C\ge L_2C$ for every irreducible
   curve $C$ in $X$.
   We assert that
   \begin{equation}\label{s rational} 
      s \not\in\bbQ \ . 
   \end{equation}
   Suppose to the contrary that $s$ is rational and consider
   the line bundle 
   $$
      L=sL_1-L_2 \in\Pic_{\bbQ}(X) \ .
   $$  
   We choose an
   integer $n$ such that $ns\in\bbZ$. The line bundle $nL$ is
   then
   algebraically equivalent to an effective (and integral)
   line bundle.  But $L$, and hence $nL$, is certainly not
   ample, so that the kernel $K(nL)$ of $\phi_{nL}$ is of
   positive dimension.  
   On the other hand, since $L_1$ and $L_2$ are not
   proportional, $nL$ is not algebraically equivalent to $0$, and
   hence $K(nL)$ cannot be the whole of $X$.  So we find that
   the neutral component of $K(nL)$ is a non-trivial abelian
   subvariety of $X$, contradicting the simplicity assumption
   on $X$.  This establishes the assertion \eqnref{s rational}

   One checks next that, since $\NX$ is finitely generated,
   its intersection with $\bbZ[L_1]\oplus\bbZ[L_2]$ is
   finitely generated as well (cf.\ for example
   \cite[Sect.\ 1.3]{Zie95}).  Choose generators
   $\liste N1k$ for this intersection.  Let now
   $0<\epsilon\ll 1$ and fix large integers $p_1,p_2$ such
   that
   \begin{equation}\label{s inequality}
      s < \frac{p_1}{p_2} < s+\epsilon \ .  
   \end{equation} 
   The line bundle 
   $$
      A=p_1L_1-p_2L_2=p_2L+p_2\(\frac{p_1}{p_2}-s\)L_1
   $$
   is then ample and therefore effective, so that we
   have $A\equiv\sum_{i=1}^k \ell_i N_i$ with integers
   $\ell_i\ge 0$.  Thus, writing $N_i\equiv a_iL_1-b_iL_2$ with
   $a_i,b_i\ge 0$, we get
   $$
      \frac{p_1}{p_2} = \frac{\sum_{i=1}^k \ell_i
      a_i}{\sum_{i=1}^k \ell_i b_i} \ ,
   $$
   which, upon letting
   $$
      q \eqdef\min\left\{\frac{a_i}{b_i}\Bigwith 1\le i \le
      k\right\} \ ,
   $$ 
   yields the lower bound
   $$
      \frac{p_1}{p_2} \ge q \ .
   $$
   But, due to the fact that 
   $s$ is irrational, which implies $q>s$, and since
   $\epsilon$ can be taken arbitrarily small, this is
   incompatible with \eqnref{s inequality}
\qed


\section{Classes on products}

   We study in this section effective divisor classes on the self-product
   $X\times X$ of an abelian variety $X$.  Suppose for a moment
   that $X$ is an elliptic curve.  Then, since $\NS(X\times X)$ is of
   rank $\ge 3$, statement (2) of the introduction says that
   $\NP(X\times X)$ is not finitely generated.  
   The argument given in \cite{Ros81} revolves around
   the alternating matrices associated
   with effective line bundles.
   Our aim here is to prove by different methods that the
   analogous statement holds in any dimension.  To simplify
   the proof, we only consider abelian varieties of Picard
   number $1$ for now, as the general case will follow with no
   effort from the proof of the theorem in Sect.\ \ref{sect proof}.

\begin{proposition}\label{product case}
   Let $X$ be an abelian variety with $\NS(X)\isom\bbZ$.  Then
   $\NP(X\times X)$ is not finitely generated.
\end{proposition}

\proof
   We denote by
   $\iota_1,\iota_2,\iota_3$ the closed embeddings of $X$ in 
   $X\times X$ given by
   $$
      \iota_1:x\mapsto(x,0),\ 
      \iota_2:x\mapsto(0,x),\
      \iota_3:x\mapsto(x,x) \ .
   $$
   Further, fix an ample line bundle $M$ whose algebraic
   equivalence class
   generates $\NS(X)$ and let
   $n$ denote the dimension of $X$.
   Supposing to the contrary that $\NP(X\times X)$ is finitely
   generated,  
   our first claim is then the following
   boundedness statement:

   \begin{items}
   \item[(*)]
      There is an integer $c>0$ such that for all effective
      line bundles $B$ on $X\times X$ with $\i B\equiv M$
      the inequality
      $$
	 \(\ii B-\iii B\)^n \le c 
      $$
      holds.
   \end{items}
   To prove (*), choose a finite set of generators $\liste N1k$
   of $\NP(X\times X)$ and write
   $$
      B\equiv\sum_{i=1}^k b_iN_i
   $$
   with integers $b_i\ge 0$.  Because of $\NS(X)=\bbZ\cdot[M]$
   we have $\i N_i\equiv n_i M$ with integers $n_i\ge 0$ for $1\le
   i\le k$.  The equivalences
   $$
      M\equiv\i B\equiv\sum_{i=1}^k b_i \i N_i
      \equiv\(\sum_{i=1}^k b_in_i\)M
   $$
   show that there is a subscript $i_0$ with the property
   $$
      b_i n_i = 
      \left\{
	 \begin{array}{ll}
	    1 & \mbox{, if $i=i_0$} \\
	    0 & \mbox{, if $i\ne i_0$} \ .
         \end{array}
      \right.
   $$
   If now $N$ is any effective line bundle on $X\times X$ with
   $\i N\equiv 0$, then it follows (for instance using the
   Seesaw Principle) that $N$ is a multiple of $pr_2^* M$,
   where $pr_2:X\times X\to X$ is the second projection, so
   that $\ii N\equiv\iii N$.  In particular we therefore obtain
   $$
      \ii B-\iii B\equiv \sum_{i=1}^k b_i\( \ii N_i-\iii N_i\)
      \equiv \ii N_{i_0} - \iii N_{i_0} \ ,
   $$
   so that (*) will hold, if we take the integer constant $c$
   to be
   $$
      c \eqdef\max\Big\{ \(\ii N_i-\iii N_i\)^n\Bigwith 
      1\le i\le k\Big\}
      \ .
   $$

   Having established (*), the idea is now to construct a
   contradiction by exhibiting a sequence of nef line bundles
   $B_m$, $m\ge 1$, satisfying
   \begin{equation}\label{B conditions}
      \i B_m\equiv M \and \lim_{m\rightarrow\infty}\(\ii
      B_m-\iii B_m \)^n =\infty \ .
   \end{equation}
   To this end we set
   $$
      L_1=pr_1^*M\ ,\ 
      L_2=pr_2^*M\ ,\
      L_3=\mu^*M\ ,
   $$
   where $pr_1,pr_2$ are the projections and $\mu$ is the
   addition map $X\times X\to X$.  We then consider the line
   bundles
   $$
      B_m\eqdef (1-m)L_1 + (m^2-m)L_2 + mL_3 \ .
   $$
   One checks that with this choice of the bundles $B_m$ the
   conditions \eqnref{B conditions}
   are satisfied.  So we will be done as soon
   as have shown that $B_m$ is nef for $m\ge 1$.

   Now recall that an ample line bundle $A$ on $X\times X$
   defines an injective homomorphism of vector spaces
   $$ 
      \NSQ(X\times X)\to\EndQ(X\times X)\ , \  
      L\mapsto\phi_A\inverse\phi_L \ ,
   $$ 
   whose image consists of the elements of $\EndQ(X\times X)$
   which are symmetric with respect to the Rosati involution
   $f\mapsto f'=\phi_A\inverse\hat f\phi_A$.  In particular,
   for an endomorphism $f$ of $X\times X$, the pullback $f^*A$
   corresponds to the symmetric endomorphism $f'f$.  Let
   now $A=L_1+L_2$.  One checks that, thanks to the fact that
   $A$ is a product polarization, an endomorphism
   $$
      f=\(\begin{array}{cc} f_1 & f_2 \\ f_3 & f_4 \end{array}\):
      X\times X\to X\times X
   $$
   is symmetric if and only if both $f_1$ and $f_4$ are
   symmetric and $f_2'=f_3$.  Therefore
   the endomorphisms
   $\alpha_1,\alpha_2,\alpha_3$, which are defined by
   \be
      \alpha_1:(x,y)&\mapsto&(x,0) \\
      \alpha_2:(x,y)&\mapsto&(0,x) \\
      \alpha_3:(x,y)&\mapsto&(x+y,x+y) 
   \ee
   are symmetric and, upon using $\alpha_1^2=\alpha_1$,
   $\alpha_2^2=\alpha_2$ and $\alpha_3^2=2\alpha_3$, 
   one finds that they correspond to the line bundles
   $L_1,L_2,L_3$.  This in turn shows that 
   the line bundle $B_m$ corresponds to the endomorphism
   $$
      \beta_m=(1-m)\alpha_1 + (m^2-m)\alpha_2 + m\alpha_3 \ .
   $$
   The point is now that $\beta_m^2=(m^2+1)\beta_m$,
   so that $\beta_m|\im\beta_m$ is just multiplication by
   $m^2+1$.  Therefore, if we denote by $Y_m$ the complementary
   abelian subvariety of $\im\beta_m$, then the differential
   of $\beta_m$ at the point $0$ is the map
   $$
      d_0\beta_m:
      T_0\im\beta_m\oplus T_0 Y_m
      \to
      T_0\im\beta_m\oplus T_0 Y_m, \
      (u,v)\mapsto((m^2+1)u,0) \ ,
   $$
   so that the analytic characteristic polynomial of $\beta_m$
   is
   $$
      P_m(t)=t^n\(t-(m^2+1)\)^n \ ,
   $$
   But the alternating coefficients of $P_m(t)$ are positive
   multiples of the intersection numbers $A^iB_m^{2n-i}$, so
   that $B_m$ is nef, as claimed.  This completes the proof of
   the proposition.
\qed
      

\section{The cone of curves and the nef cone
   of an abelian variety}\label{sect proof}

   Finally, we give in this section the proof of
   the theorem stated in the introduction.
   So let $X$ be an abelian variety and 
   denote by $N_1(X)$ the
   vector space of numerical equivalence classes of
   real-valued 1-cycles on $X$, and by $\NE(X)$ the convex
   cone in $N_1(X)$ generated by irreducible curves.  Through
   the intersection product the vector space $N_1(X)$ is dual
   to the N\'eron-Severi vector space
   $\NSR(X)=\NS(X)\tensor\bbR$.  The dual cone of $\NE(X)$ is
   the nef cone 
   $$
      \Nef(X)
      =\{\lambda\in\NSR(X)\with\lambda\xi\ge 0\mbox{
      for all }\xi\in\NE(X)\} \ ,
   $$
   which in the case of abelian varieties coincides with the
   effective cone (cf.\ Lemma \ref{nef and effective}).
   The dual of $\Nef(X)$ in turn is the closed cone $\bNE(X)$,
   so that
   one of these two cones is rational polyhedral if and only if the
   other is.  By Gordon's Lemma this is 
   equivalent to the semi-group $\NX$ being finitely
   generated.  (See e.g.\ \cite[Theorem 14.1 and \S\S19,20]{Roc70} 
   for the elementary properties of cones used here).

   The idea is now, given an abelian variety, to first apply 
   Poincar\'e`s Complete Reducibility
   Theorem, i.e.\ to decompose it up to
   isogenies into a product of powers of non-isogenous simple abelian
   varieties, and to apply Proposition \ref{simple case} and
   Proposition \ref{product case} subsequently.  
   One needs then that 
   finite generation of $\NX$ is a property which is invariant
   under isogenies:

\begin{lemma}[\cite{Ros81}]\label{isog}
   Let $X$ and $Y$ be isogenous abelian varieties.  Then
   $\NP(X)$ is finitely generated if and only if $\NP(Y)$ is.
\end{lemma}

   Since this observation is crucial for our approach, 
   let us briefly indicate a proof, before
   we proceed to the proof of the theorem.  
   So suppose that $\NX$ is finitely generated and that there
   is an isogeny $f:X\to Y$.  Thanks to the fact that $f^*$
   embeds $\NP(Y)$ into $\NX$ and to the symmetry of the
   situation, it is enough to show that $f^*\NP(Y)$ is finitely
   generated.  Let then $\liste N1k$ be generators for $\NX$
   and put for $1\le i\le k$
   $$
      n_i \eqdef\min\{ n\in\bbZ\with nN_i\in f^*\NP(Y) \} \ .
   $$
   (The set on the right-hand side is non-empty, since $f$ is
   an isogeny.)  Then $f^*\NP(Y)$ is generated by the elements
   $n_1N_1,\dots,n_rN_r$ together with those elements
   $\sum_{i=1}^k m_iN_i$, $0\le m_i<n_i$, which belong to
   $f^*\NP(Y)$.
\qed

   Consider now for an abelian variety $X$ its decomposition 
   $$
      X_1^{n_1}\times\dots\times X_r^{n_r}
   $$
   up to isogenies, where $\liste X1r$
   are mutually non-isogenous
   simple abelian varieties and $\liste n1r$ are positive
   integers.  In view of the remarks made at the beginning of
   this section, the theorem stated in the introduction will 
   follow from

\begin{theorem}
   The semi-group $\NX$ is finitely generated if and only if
   $\NS(X_i)\isom\bbZ$ and $n_i=1$ for $1\le i\le r$.
\end{theorem}

\proof
   Suppose first that the conditions on the factors $X_i$ and the
   exponents $n_i$ are satisfied for $X$.  
   By Lemma \ref{isog} we may assume that $X$
   is the product $X_1\times\dots\times X_r$.  
   Fix for $1\le i\le r$ an ample 
   generator $N_i$ of $\Pic(X_i)$ and let $A$ be the
   product polarization $A=\sum_{i=1}^r pr_i^*N_i$.  Due to
   the fact that the $X_i$ are non-isogenous, we have
   $$
      \NSQ\(\prod_{i=1}^r X_i\)
      \isom\EndQ^s\(\prod_{i=1}^r X_i\)
      \isom\bigoplus_{i=1}^r\EndQ^s(X_i)
      \isom\bigoplus_{i=1}^r\NSQ(X_i) \ ,
   $$
   where $\EndQ^s\(\prod_{i=1}^r X_i\)$ and $\EndQ^s(X_i)$
   denote the subgroups of symmetric endomorphisms with
   respect to the Rosati involutions associated with $A$ and
   $N_i$ respectively.  Therefore
   $$
      \NX=\bigoplus_{i=1}^r \bbZ^+\cdot[N_i] 
   $$ 
   is finitely generated.

   Now suppose conversely 
   that $\NX$ is finitely generated.  
   By Lemma \ref{isog} again, we may assume that
   $X$ is the product $X_1^{n_1}\times\dots\times X_r^{n_r}$.
   Note that if $V_1$ and $V_2$ are
   varieties such that $\NP(V_1\times V_2)$ is finitely generated,
   then $\NP(V_1)$ and $\NP(V_2)$ are finitely generated as
   well.  So in particular Proposition 
   \ref{simple case} applies to the factors $X_i$ and shows
   that we have $\NS(X_i)\isom\bbZ$ for all $i$.  
   Further, if we had $n_i>1$ for some $i$, i.e.\ if
   a multiple factor $X_i$ appeared in the product decomposition of
   $X$, then $\NP(X_i\times X_i)$
   would be finitely generated, which however is impossible
   according to Proposition \ref{product case}.  
   This completes the proof
   of the theorem.
\qed

\begin{varthm*}{\it Acknowledgements}
\rm
   This research was supported by DFG contract Ba 1559/2-1.
   I would like to thank R.\ Lazarsfeld for helpful
   discussions and the University of California, Los Angeles,
   for its hospitality.
\end{varthm*}


\bigskip
\bigskip

   Thomas Bauer,
   Department of Mathematics,
   University of California, 
   Los Angeles, CA 90095--1555

   (E-mail: {\tt tbauer@math.ucla.edu})

   {\em Current address:}
   Mathematisches Institut,
   Universit\"at Erlangen-N\"urnberg,
   Bis\-marck\-stra{\ss}e $1\frac12$,
   D-91054 Erlangen,
   Germany

   (E-mail: {\tt bauerth@mi.uni-erlangen.de})


\end{document}